\documentclass[12pt]{article}
\usepackage[cp1251]{inputenc}
\usepackage[russian]{babel}
\usepackage{amsmath,amsthm,amssymb}
\usepackage{latexsym}
\usepackage{amsfonts}
\newcommand{\en}{\enspace}

\begin{document}

$$
$$

\begin{center}
\bf SCATTERING AND BOUND STATES FOR\\
NONSELFADJOINT SCHR\"ODINGER OPERATOR
\end{center}

\medskip

\begin{center}
\bf S. A. Stepin
\end{center}

\bigskip
\bigskip

\textbf{Abstract.}
 Spectral components of one-dimensional Schr\"odinger operator with complex potential are investigated. An effective upper bound for the total number of eigenvalues and spectral singularities is established. For dissipative Schr\"odinger operator sufficient condition is found which guarantees the absence of singular component in the continuous spectrum and spectral decomposition is exposed.

\medskip


\medskip


\bigskip
\bigskip

\begin{center}
\bf \S 1. Introduction
\end{center}

\medskip

Characteristic features which distinguish nonselfadjoint operators with continuous spectrum from selfadjoint ones are exhibited in the case of one-dimensional Schr\"odinger operator $\,L=-d^2/dx^2\,+\,V(x)\,$ on half-line with complex potential $\,V(x)\,$ and Dirichlet boundary condition at zero. Such an operator proves to be quite a simple and rather rich model which displays a number of effects typical for perturbation theory in nonselfadjoint setting (see [1] and [2]). Besides that, Schr\"odinger operator with complex potential is known to appear (see [3]) in the study of open quantum mechanical systems with energy dissipation.

The problem about similarity of the part $\,L_c\,$ of operator $\,L\,$ associated with its continuous spectrum and an unperturbed operator $\,L_0\,$ corresponding to $\,V(x)\equiv 0\,$ is closely related with generalized eigenfunction expansion problem.
In the nonselfadjoint case under the condition of local integrability of the potential $\,V(x)\,$ the corresponding spectral function is sentenced in [4] to be a singular distribution.
When imaginary part of $\,V(x)\,$ is bounded such a generalized spectral function can be represented as a contour (containing the spectrum of the problem inside) integral; sometimes this contour admits a deformation right to the spectrum itself which gives an eigenfunction expansion of classical form (see [5]).

For Sturm-Liouville problems on half-line the so-called transformation operators prove to be effective tools in the study of spectral similarity as well as for the solving of direct and inverse scattering problems. In turn scattering theory itself provides an adequate construction of transformation operators intertwining perturbed operator with an unperturbed one. This approach elaborated in [6] and [7] for Schr\"odinger operator with complex potential is based on nonstationary definition of wave operators realizing similarity of $\,L_c\,$ and $\,L_0.\,$ To this end wave operators are constructed by means of comparison of the corresponding perturbed and unperturbed dynamics and characterize the first one at large time periods. Existence of complete wave operators implies the absence of singular component in the continuous spectrum of $\,L\,$ in the sense that $\,L_c\,$ is similar to an operator with absolutely continuous spectrum.

Certain difficulties in the study of general nonselfadjoint operators and Schr\"odinger operator with complex potential in particular are caused by the lack of a priori information about the behavior of the resolvent near the points of spectrum as well as complicated structure of the spectrum itself. In this context an essential role is played by the spectral singularities, i.e. the poles of the analytic continuation of the resolvent integral kernel which are embedded into continuous spectrum. It turns out that the key information about spectral properties of operator $\,L\,$ can be extracted from the so-called Jost function which coincides under appropriate assumptions with Fredholm determinant
$$
e(k)\en=\en\det\big(\,I\,+\,V\,|V|^{-1/2}(L_0-k^2I)^{-1}|V|^{1/2}\big)\,;
$$
for example Nevanlinna-Riesz canonical factorization of Jost function generates a decomposition of the underlying Hilbert space into the sum of invariant (with respect to operator $\,L)\,$ subspaces which correspond to different components of the spectrum.

In the present paper the spectral components of one-dimensional Schr\"odinger operator with complex potential are investigated. Theorem 1 gives an effective upper bound for the total number of eigenvalues and spectral singularities of operator $\,L\,$ under certain assumptions imposed on the potential which allow the resolvent integral kernel to possess an analytic continuation through the continuous spectrum cut. Such estimates of integral type have been known previously in the selfadjoint case only when a spectral singularity can appear on the continuous spectrum edge and corresponds to the so-called virtual level.

Provided that potential $\,V(x)\,$ possesses the first momentum and its imaginary part is positively definite a sufficient condition for the corresponding dissipative Schr\"odinger operator is obtained which guarantees the absence of singular component in the spectrum and by the usage of wave operators a spectral decomposition for $\,L\,$ is given. The class of dissipative operators was introduced and studied in an abstract setting by R.Phillips [8] and K.Friedrichs [9]; as regards Schr\"odinger operator the importance of this class in the context of the eigenfunction expansion problem was manifested in [1] (see also [7]).

Inversion formulas derived in [2] for the generalized Fourier-Laplace transform generating spectral representation of $\,L\,$ enable one to construct transformation operators realizing similarity of $\,L_c\,$ and $\,L_0\,$ however under additional condition imposed on $\,V(x)\,$ involving finiteness of its second momentum. Theorem 2 below can be treated as an extension (and a progress in certain aspects) of the results established in [2], [6] and [7].
The sufficient conditions for absolute continuity of the spectrum of operator $\,L_c\,$ enlarge the class of potentials $\,V(x)\,$ for which a generalized spectral function characterized in [4] turns out to be regular while a construction of wave operators intertwining $\,L_c\,$ and $\,L_0\,$ enables one to obtain a useful information about the structure of the corresponding spectral measure.

The present paper is organized as follows. In \S\,2 we define the required notions and introduce the subjects which are in the scope, further the formulations of the main results (theorems 1 and 2) are presented as well as their comparison with the previous investigations is carried out. In \S\,3 we study the discrete spectrum of operator $\,L\,$ and derive an estimate for the total number of its eigenvalues and spectral singularities (proof of theorem 1). Section 4 plays an important preparatory role. In statement 1 we establish some facts to be used about the transformation operator for the pair $\,L_0\,$ and $\,\widehat{L}={\rm Re}\,L\,$ and crucial estimates of its integral kernel are obtained (lemmas 1 and 2). In \S\,5 the description of spectral components for the dissipative operator $\,L$ is given (statement 3), based on the exploiting in the statement 2 the analytic properties of the corresponding Jost function. Section 6 deals with the construction of direct and inverse wave operators for the pair $\,\widehat{L}\,$ and $\,L$ (statement 4). The key role here is played by the property of smooth vectors to be dense in the continuous spectrum subspace of operator $\,\widehat{L}\,$ which we establish in lemma 4 making use of the results of \S\,4. Section 7 gathers and sums up the previous considerations. Therein we finish the proof of theorem 2 and obtain a simple sufficient condition for the operator $\,L\,$ to be completely nonselfadjoint. In conclusion stationary representations of wave operators are indicated (statement 5) which clarify the relationship with eigenfunction expansion for operator $\,L.$

\bigskip

\begin{center}
\bf \S\,2. Statement of main results
\end{center}

\bigskip

In the space $\,{\cal H}={\rm L}_2(0,\infty)\,$ consider an operator
$$
L\en=\en L_0+V\en=\en -\,d^2/dx^2\,+\,V(x)
$$
with the domain $\,D(L)\,=\,\{\,y\in{\cal H}\!:\,y'$ {\it absolutely continuous}, $y''\in{\cal H},\,y(0)=0\,\},\,$ where potential $\,V(x)\,$ is bounded and complex-valued. When $\,V(x)\equiv 0\,$ the corresponding operator will be denoted by $\,L_0.$

Under the standard scattering theoretic condition imposed on the potential --- the finiteness of its first momentum, the absence of spectral singularities, i.e. the points of continuous spectrum corresponding to zeroes of Jost function (cf. \S\,3), guarantees that the set of eigenvalues of operator $\,L\,$ is finite (see [1]). On the other hand provided that potential $\,V(x)\,$ is integrable on half-line with exponential weight the corresponding Jost function admits analytic continuation through continuous spectrum cut and therefore eigenvalues have no accumulation points. The monograph [10] contains a detailed survey of early period papers in which various sufficient conditions can be found providing quasi-analyticity of Jost function and thus finiteness of discrete spectrum of $\,L.\,$ As regards the later period results in this direction paper [11] should be mentioned where finiteness of the set of eigenvalues is established under the assumption that potential $\,V(x)\,$ admits a bounded analytic continuation into a sector containing half-line $\,\mathbb R_+\,$ and possesses finite first momenta on its sides.

The estimate of the total number $\,N(V)\,$ of eigenvalues and spectral singularities for operator $\,L\,$ is given by

\bigskip

{\bf Theorem 1.} {\it Provided that $\,V\in{\rm L}_1(0,\infty)\,$ all the eigenvalues and spectral singularities of operator $\,L\,$ are located in the disc $\,\,\,|\lambda\,|\,\leqslant\,R\,,\,\,$ of radius
$$
R\,=\,R(\alpha)\,:=\,\bigg(\frac1{\ln 2}\int_0^{\infty}\!\!(2x)^{\alpha}|V(x)|\,dx\bigg)^{2/(1-\alpha)}\!,\quad \alpha\in[0,1)\,.
$$
Given $\,a>0\,$ suppose that the integral $\,\,\displaystyle{\int^{\infty}\!e^{ax}|V(x)|\,dx}\,\,$ converges, then
 the inequality
$$
N(V)\,\,\leqslant\,\,\bigg(\!\ln\frac{A+a/2}{\sqrt{A^2+R}}\bigg)^{-1}
\bigg\{\,\frac1{a^{1-b}}\int_0^{\infty}\!\!x^b(1+e^{ax})|V(x)|\,dx\,-\,\ln\Big(2-2^{(\sqrt{R}/A)^{1-\alpha}}\Big)\bigg\}
$$
is valid for arbitrary $\,A>\max\{\sqrt{R},R/a-a/4\}\,$ and $\,b\in[0,1].$
}

\bigskip

Putting $\,\,c=\min\{a,\sqrt{R}\}\,\,$ and choosing appropriate values of parameters $\,A\,$ and $\,b\,$ we obtain the following estimate
$$
N(V)\,\,\leqslant\,\,\frac2{\ln\,(1+c^2/4R)}\,
\bigg\{\,\frac1{c}\int_0^{\infty}\!\!(1+e^{cx})|V(x)|\,dx\,-\,\ln\!\Big(2-2^{(c/2\sqrt{R})^{1-\alpha}}\Big)\bigg\}\,.
$$

Previously such estimates for the quantity $\,N(V)\,$ were known in selfadjoint case only when Birman-Schwinger principle can be employed to this end. The upper bound for $\,N(V)\,$ from theorem 1 is obtained by using Jensen formula applied to counting function for the number of zeroes of $\,e(k).\,$  Operators of the type considered here are known to appear in singular perturbation theory as a model describing the passage from discrete spectrum to continuous one. In [12] the distribution of eigenvalues for such operators with purely imaginary potential was investigated and within quasi-classical approach to localization of the spectrum the corresponding Bohr-Sommerfeld quantization rules were derived.

Continuous spectrum component of dissipative Schr\"odinger operator with purely imaginary potential was investigated in [13]. In the present paper the general case is considered when potential $\,V(x)\,$ along with imaginary part may have nontrivial real one as well. For the sake of simplicity we shall assume that operator $\,L\,$ is completely nonselfadjoint, i.e. does not possess any non-trivial selfadjoint  restriction.

\bigskip

{\bf Theorem 2.} {\it Suppose that completely nonselfadjoint dissipative operator $\,L=L_0+V\,$ with bounded potential $\,V(x),\,$ satisfying the condition
\begin{equation*}
\int_0^{\infty}x\,|\,V(x)|\,\,dx\en<\en\infty\,\,,
\end{equation*}
has no spectral singularities. Then the bounded direct wave operator
$$
\Omega\en=\en{\rm s}\!-\!\!\lim_{t\to\infty}\exp(itL)\,\exp(-itL_0)
$$
exists and possesses completeness property. On the invariant subspace $\,{\cal H}_c=\Omega{\cal H}\,$ corresponding to continuous spectrum of $\,L\,$ the inverse wave operator
$$
\widetilde\Omega\en=\en{\rm s}\!-\!\!\lim_{t\to\infty}\exp(itL_0)\,\exp(-itL)
$$
is well-defined and bounded, so that $\,L_c=L|{\cal H}_c\,$ and $\,L_0\,$ are similar\,{\rm :} $\,\,L_c\,=\,\Omega\,L_0\,\widetilde\Omega\,.$
}

\bigskip

In the context of theorem 2 note that a dissipative operator $\,L=L_0+V\,$ is completely nonselfadjoint if there exists an interval on which $\,{\rm Im}\,V(x)>0;\,$ this condition guarantees that discrete spectrum of $\,L\,$ is located in the open upper half-plane $\,\mathbb C_+$ (cf. [1]).

Operator $\,\widetilde\Omega\,$ proves to be the left inverse for $\,\Omega\,$ and its right inverse on the subspace $\,{\cal H}_c.$
Thus theorem 2 under appropriate assumptions gives solution to the problem (see [7]) of constructing the spectral representation for the operator in question; to this end we make use of direct and inverse wave operators $\,\Omega\,$ and $\,\widetilde\Omega\,$ to implement similarity of $\,L_0\,$ and restriction $\,L_c\,$ of operator $\,L\,$ onto the continuous spectrum subspace $\,{\cal H}_c.\,$ Note that direct wave operator $\,\Omega\,$ exists and is bounded without regard for the absence (or presence) of spectral singularities in the continuous spectrum of dissipative operator $\,L.\,$

The first momentum of the potential $\,V(x)\,$ is required to be finite in the selfadjoint case (see [14]) to construct eigenfunction expansion for operator $\,L\,$ and to derive the corresponding inversion formula. Continuous spectrum eigenfunctions of operator $\,L\,$ are naturally obtained from those of operator $\,L_0\,$ under the action of wave operator while the absolutely continuous part of the spectral measure proves to be of Lebesgue type. Given locally integrable potential $\,V(x)\,$ one can only state (cf. e.g. [15]) the very fact of existence of the corresponding spectral function which (as it was already mentioned) may occur singular in nonselfadjoint case.

Under the condition
$$
\int_0^{\infty}(1+x^2)\,|\,V(x)|\,\,dx\en<\en\infty
$$
eigenfunction expansion formulas for a nonselfadjoint operator $\,L=L_0+V\,$ derived in [1] and [2] enable one to obtain transformation operators realizing similarity of $\,L_c\,$ and $\,L_0.\,$ The required similarity was also established in [6] provided that
$$
\int_0^{\infty}x\,|\,V(x)|\,\,dx\en<\en 1\,.
$$
A dissipative Schr\"odinger operator on the half-line under the condition
$$
\int_0^{\infty}x^{1+\varepsilon}\,|V(x)|\,\,dx\en<\en\infty\,,\quad\varepsilon>0\,,
$$
(which is somewhat more restrictive than the finiteness of the first momentum) was studied in [7] where the direct wave operator for the pair $\,\{L, L_0\}\,$ was constructed and its relationship with the transformation operator mentioned above was discussed (cf. \S\,7).

The results of the present paper extend and supplement the previous investigations. A characteristic feature which distinguishes our approach is the treatment of $\,{\rm Re}\,L\,$ as an unperturbed operator for $\,L\,$ and its comparison with $\,L_0\,$ by means of Volterra type transformation operator. The estimates of its integral kernel (see \S\,4) enable us to establish for the perturbation the smoothness property of Kato type with respect to (appropriately chosen) unperturbed operator. Due to this fact as well as to the absence of singular component in Nevanlinna-Riesz canonical factorization of Jost function the completeness property for $\,\Omega\,$ and boundedness of  $\,\widetilde\Omega\,$ are valid provided that operator $\,L\,$ has no spectral singularities.

\bigskip

\begin{center}
\bf \S\,3. Discrete spectrum of operator $L$
\end{center}

\medskip

Given integrable potential $\,V(x)\,$ the equation
\begin{equation}
-\,y''\,+\,V(x)\,y\en=\en k^2y
\label{form1}
\end{equation}
is known (see [1]) to possess for any $\,k\in{\mathbb C_+}\,$ a solution $\,e(x,k)\,$ with asymptotics
$$
e(x,k)\en\sim\en e^{ikx}\,,\quad x\to\infty\,.
$$
Furthermore the spectrum of operator $\,L\,$ consists of continuous and discrete components
$$
\sigma_c(L)\,=\,\mathbb R_+\,,\quad \sigma_d(L)\,=\,\{\,k^2\!:\,e(k)=0,\,k\in\mathbb C_+\,\}\,,
$$
where $\,e(k):=e(0,k)\,$ is the so-called Jost function. The set of eigenvalues $\,\sigma_d(L)\,$ is bounded, at most countable and its accumulation points (if any) belong to the half-line $\,\mathbb R_+,\,$ whereas operator $\,L\,$ has no eigenvalues embedded into continuous spectrum.

Under the condition
\begin{equation}
\int_0^{\infty}\!\! x|V(x)|\,dx\en<\en\infty
\label{form2}
\end{equation}
function $\,e(k)\,$ is analytic in the open half-plane $\,\mathbb C_+\,$ and continuous up to its boundary.
Real zeroes of Jost function $\,e(k)\,$ correspond to distinguished points $\,\lambda=k^2\,$ of continuous spectrum $\,\sigma_c(L)\,$ which are known as {\it spectral singularities}. Provided that operator $\,L\,$ has no spectral singularities its discrete spectrum clearly proves to be  finite.
An upper bound for the total number $\,N(V)\,$ of eigenvalues and spectral singularities of operator $\,L=L_0+V\,$ will be obtained here under additional restriction
\begin{equation}
\int^{\infty}e^{ax}|V(x)|\,dx\,\,<\,\,\infty\,.
\label{form3}
\end{equation}

\medskip

{\bf Proof of theorem 1.} Under the condition (\ref{form3}) function $\,e(k)\,$ admits analytic continuation from $\,\mathbb C_+\,$ to the strip $\,\Pi(a)=\{-a/2<{\rm Im}k\leqslant 0\}.$ Indeed Jost solution satisfies the integral equation
$$
e(x,k)\en=\en e^{ikx}\,-\,\int_x^{\infty}\frac{\sin k(x-\xi)}{k}\,V(\xi)\,e(\xi,k)\,d\xi
$$
and consequently is represented by the series $\,\displaystyle{e(x,k)=e^{ikx}\sum_{n=0}^{\infty}\varepsilon^{(n)}(x,k)}\,$
where $\,\varepsilon^{(0)}(x,k)=1\,$ and
$$
\varepsilon^{(n)}(x,k)\en=\en\int_x^{\infty}\frac{e^{2ik(\xi-x)}-1}{2ik}\,V(\xi)\,\varepsilon^{(n-1)}(\xi,k)\,d\xi\,.
$$
By induction in $\,n\,$ and with the usage of inequality $\,\,\big|\,e^{2ikx}-1\big|\,\leqslant\,(2|k|x)^{\alpha}\big(1+e^{-2{\rm Im}k\,x}\big),$ $\,\alpha\in[0,1],\,$
one can verify for $\,k\in\mathbb C_+\!\cup\Pi(a)\,$ the following estimate
$$
|\,\varepsilon^{(n)}(x,k)\,|\,\,\leqslant\,\,\frac1{n!}\,\bigg(\frac1{(2|k|)^{1-\alpha}}\int_x^{\infty}\!\!\xi^{\alpha}\big(1+e^{-2{\rm Im}k\,\xi}\big)|V(\xi)|\,d\xi\bigg)^n.
$$
Thus the above series for $\,e(k)=e(0,k)\,$ converges uniformly in $\,\mathbb C_+\!\cup\Pi(a),\,$ therefore it gives an analytic function and besides for arbitrary $\,\alpha\in[0,1]\,$ the inequality
\begin{equation}
\big|\,e(k)-1\,\big|\,\,\leqslant\,\,\exp\bigg(\frac1{(2|k|)^{1-\alpha}}\int_0^{\infty}\!\!\xi^{\alpha}\big(1+e^{-2{\rm Im}k\,\xi}\big)|V(\xi)|\,d\xi\bigg)-\,1\,\,
\label{form4}
\end{equation}
holds.
Note that Jost function admits (see [16]) a representation in the form of Fredholm determinant
$$
e(k)\,\,=\,\,\det\Big(\,I\,+\,V\,|V|^{-1/2}\big(L_0-k^2I\big)^{-1}|V|^{1/2}\Big)\,,
$$
in which the resolvent $\,(L_0-k^2I)^{-1}$ is an integral operator with the kernel $\,\exp\big(ik\max\{x,\xi\}\big)\sin\big(k\min\{x,\xi\}\big)/k\,,\,$
where the required analytic continuation of $\,e(k)\,$ clearly results from.

\smallskip

The total number $\,N(V)\,$ of eigenvalues and spectral singularities of operator $\,L\,$ coincides with the multiplicity of the zeroes of Jost function $\,e(k)\,$ in the closed upper half-plane $\,\overline{\mathbb C}_+,\,$ while by virtue of (\ref{form4}) all of them are located in the domain
$$
\big\{k\in\overline{\mathbb C}_+:|\,k|\leqslant r\big\}\,\,\subset\,\,\big\{|\,k-iA|\leqslant\sqrt{A^2+r^2}\,\big\}\,,
$$
depending on the potential $\,V\,$ by means of parameters
$$
r\,=\,r(\alpha)\,=\,\bigg(\frac1{\ln 2}\int_0^{\infty}\!\!(2x)^{\alpha}|V(x)|\,dx\bigg)^{1/(1-\alpha)}\!\!,\quad
A>\max\{r,r^2/a-a/4\}\,.
$$
To evaluate the quantity $\,N(V)\,$ which does not exceed the number of zeroes of the function $\,\varphi(z):=e(z+iA)\,$ located in the disc $\,|z|\leqslant\sqrt{A^2+r^2}\,$ we make use of Jensen formula. Beforehand note that $\,\sqrt{A^2+r^2}<A+a/2\,$ since $\,A>r^2/a-a/4.$ Choosing an arbitrary $\,R\in(\sqrt{A^2+r^2},A+a/2)\,$ and setting
$$
n_{\varphi}(t)\,\,:=\,\,\#\,\big\{z_j:\,\varphi(z_j)=0, |z_j|<t\big\}
$$
to denote zeroes counting function for $\,\varphi\,$ we obtain
\begin{multline*}
N(V)\,\ln\frac{R}{\sqrt{A^2+r^2}}\,\,\leqslant\sum_{|z_j|^2\leqslant\,A^2+r^2}\ln\frac{R}{|z_j|}\,\,=\,\,\int_0^R\frac{n_{\varphi}(t)}{t}\,dt\,\,=\\
=\,\,\frac1{2\pi}\int_0^{2\pi}\!\ln|\varphi(R e^{i\theta})|\,d\theta\,\,-\,\,\ln|\varphi(0)|\,.
\end{multline*}

\smallskip

In order to estimate the absolute value of $\,\varphi\,$ on the circle $\,|z|=R\,$ inequality (\ref{form4}) is applied\,:
\begin{multline*}
|\varphi(R e^{i\theta})|\,=\,|\,e(iA+R e^{i\theta})|\,\,\leqslant\,\,\exp\bigg(\frac1{(2|\,iA+R e^{i\theta}|\,)^{1-b}}\int_0^{\infty}\!\xi^b\Big(\,1+e^{-2(A+R\sin\theta)\xi}\Big)|V(\xi)|\,d\xi\bigg)\\
\leqslant\,\,\exp\bigg(\frac1{(2(R-A))^{1-b}}\int_0^{\infty}\!\xi^b\Big(\,1+e^{2(R-A)\xi}\Big)|V(\xi)|\,d\xi\bigg)
\end{multline*}
where $\,b\in[0,1].\,$ Again by virtue of (\ref{form4}) for $\,\alpha\in[0,1]\,$ one has
\begin{multline*}
|\varphi(0)|\,\geqslant\,1\,-\,|e(iA)-1|\,\geqslant\,2\,-\,\exp\bigg[\frac1{(2A)^{1-\alpha}}\int_0^{\infty}\!x^{\alpha}(1+e^{-2Ax})|V(x)|\,dx\bigg]\,\,\geqslant\\
\geqslant\,\,2\,-\,\exp\bigg[\frac1{A^{1-\alpha}}\int_0^{\infty}\!(2x)^{\alpha}|V(x)|\,dx\bigg]\,\,=\,\,2-2^{(r/A)^{1-\alpha}}\,,
\end{multline*}
where the right-hand-side is positive due to the choice $\,A>r.$
Thus for arbitrary $\,R\in(\sqrt{A^2+r^2},A+a/2)\,$ an inequality
\begin{multline*}
N(V)\,\ln\frac{R}{\sqrt{A^2+r^2}}\,\,\leqslant\,\,
\frac1{(2(R-A))^{1-b}}\int_0^{\infty}\!\xi^b\Big(1+e^{2(R-A)\xi}\Big)|V(\xi)|\,d\xi\,\,\,-\\
-\,\,\ln\Big(2-2^{(r/A)^{1-\alpha}}\Big)
\end{multline*}
is proved to be true. Passing here to the limit as $\,R\to A+a/2\,$ we get the desired estimate for the number of eigenvalues and spectral singularities of operator $\,L:$
$$
N(V)\,\ln\frac{A+a/2}{\sqrt{A^2+r^2}}\,\,\leqslant\,\,
\frac1{a^{1-b}}\int_0^{\infty}\!\!\xi^b(1+e^{a\xi})|V(\xi)|\,d\xi\,-\,\ln\Big(2-2^{(r/A)^{1-\alpha}}\Big)\,.
$$

\smallskip

Setting $\,\mu=\min\{a,r\}\,$ and choosing $\,b=0,\,$ for $\,A>\max\{r,r^2/\mu-\mu/4\}\,$ we obtain the following estimate
$$
N(V)\,\ln\frac{A+\mu/2}{\sqrt{A^2+r^2}}\,\,\leqslant\,\,
\frac1{\mu}\int_0^{\infty}\!\!(1+e^{\mu\xi})|V(\xi)|\,d\xi\,-\,\ln\Big(2-2^{(r/A)^{1-\alpha}}\Big)\,.
$$
In particular when $\,A=2r^2/\mu\,$ one has
$$
N(V)\,\,\leqslant\,\,2\,\Big(\!\ln\Big(1+\mu^2/4r^2\Big)\Big)^{-1}
\bigg\{\,\frac1{\mu}\int_0^{\infty}\!\!(1+e^{\mu x})|V(x)|\,dx\,-\,\ln\Big(2-2^{(\mu/2r)^{1-\alpha}}\Big)\bigg\}\,.
$$

\bigskip

{\bf Corollary 1.} \, {\it
If $\,\,\displaystyle{a\,\geqslant\, \rho\,:=\,\frac1{\ln 2}\int_0^{\infty}\!\!|V(x)|\,dx}\,\,$ then the estimate holds
\begin{multline*}
N(V)\en\leqslant\en\frac2{\ln(5/4)}\,\bigg\{\,\frac1{\rho}\int_0^{\infty}(1+e^{\rho x})|V(x)|\,dx\,+\,\ln\bigg(1+\frac1{\sqrt{2}}\bigg)\bigg\}\en\leqslant\\
\leqslant\en 10\,\bigg\{\,1\,+\,\frac2{\rho}\int_0^{\infty}e^{\rho x}|V(x)|\,dx\bigg\}\,.
\end{multline*}
}

\bigskip

\begin{center}
\bf \S\,4. Transformation operator for the pair $\,L_0\,$ and $\,{\rm Re}\,L$
\end{center}

\medskip

Let us put into usage a notation $\,p(x)\,:=\,{\rm Re}\,V(x)\,$ and along with operator $\,L=L_0+V=-d^2/dx^2+V(x)\,$ generated in $\,{\cal H}={\rm L}_2(0,\infty)\,$ by the boundary condition $\,y(0)=0\,$ consider a selfadjoint operator $\,\widehat{L}=-d^2/dx^2+p(x)\,$ with the same condition at zero. If $\,p\in{\rm L}_1(0,\infty)\,$ then $\,\sigma_c(\widehat{L})=\mathbb R_+\,$ while operator $\,\widehat{L}\,$ has no embedded eigenvalues; besides, the condition
\begin{equation}
\int^{\infty}\!x|p(x)|\,dx\,<\,\infty
\label{form5}
\end{equation}
ensures that also the edge of continuous spectrum is not an eigenvalue of operator $\,\widehat{L}.\,$ Let $\,\widehat{\cal H}_c\,$ be an absolutely continuous spectrum subspace of operator $\,\widehat{L}\,$ and set $\,\widehat{L}_c:\,=\,\widehat{L}\,|\,\widehat{\cal H}_c\,.$

For $\,k\in\mathbb C_+\,$ denote by $\,s(x,k)\,$ a solution to equation
\begin{equation}
-\,y''\,+\,p(x)\,y\en=\en k^2y
\label{form6}
\end{equation}
satisfying the boundary conditions $\,\,s(0,k)=0\,,\,\, s'_x(0,k)=k.\,$
Solution $\,s(x,k)\,$ admits a representation
$$
s(x,k)\en=\en\sin kx\en+\en\int_0^xK(x,\xi)\,\sin k\xi\,d\xi
$$
so that the following assertion is valid.

\bigskip

{\bf Lemma 1.} {\it The integral kernel $\,K(x,\xi)\,$ satisfies the estimate
$$
|K(x,\xi)|\en\leqslant\en\frac12\int_{(x-\xi)/2}^x|p(r)|\,dr\,\exp\bigg(\frac12\int_0^{x-\xi}\!\!\bigg(\int_{s/2}^{\infty}\!|p(r)|\,dr\bigg)ds\bigg)\,.
$$
}

\smallskip

{\bf Proof.} Kernel $\,K(x,\xi)\,$ proves to be (see e.g. [17]) a solution to the integral equation
$$
K(x,\xi)\en=\en\frac12\int_{(x-\xi)/2}^{(x+\xi)/2}p(u)\,du\,\en+\en\int_{(x-\xi)/2}^{(x+\xi)/2}
du\int_0^{(x-\xi)/2}\!\!\!p(u+v)\,K(u+v,u-v)\,dv\,,
$$
which takes in the new variables $\,\,z=x+\xi,\,\,w=x-\xi\,\,$ the form
$$
A(z,w)\en=\en\frac12\int_{w/2}^{z/2}p(t)\,dt\en+\en\frac14\int_w^zdt\int_0^wp((t+s)/2)\,A(t,s)\,ds
$$
where $\,\,A(z,w)\,=\,K((z+w)/2,(z-w)/2).\,$ Solving this equation by iterations we get
\begin{eqnarray*}
A(z,w)&=&\sum_{n=0}^{\infty}A_n(z,w)\,,\en A_0(z,w)\,=\,\frac12\int_{w/2}^{z/2}p(t)\,dt\,,\\
A_n(z,w)&=&\frac14\int_w^zdt\int_0^wp((t+s)/2)\,A_{n-1}(t,s)\,ds\,.
\end{eqnarray*}
By induction argument the following estimate
$$
|A_n(z,w)|\en\leqslant\en\frac1{2n!}\int_{w/2}^{(z+w)/2}\!\!|p(r)|\,dr\,\,
\bigg(\frac12\int_0^w\bigg(\int_{s/2}^{(z+s)/2}\!\!|p(r)|\,dr\bigg)\,ds\bigg)^n
$$
is then established. Indeed changing the order of integration and substituting $\,t+s=2r\,$ one obtains
\begin{multline*}
|A_{n+1}(z,w)|\en\leqslant\en\frac14\int_0^wds\int_w^z|p((t+s)/2)|\,|A_n(t,s)|\,dt\en\leqslant\\
\leqslant\en\frac1{8n!}\int_0^wds\int_w^z|p((t+s)/2)|\,dt\int_{s/2}^{(z+s)/2}\!\!|p(r)|\,dr\,\,\bigg(\frac12\int_0^s\bigg(\int_{r/2}^{(z+r)/2}\!\!|p(u)|\,du\bigg)\,dr\bigg)^n\,\leqslant\\
\leqslant\en\frac1{4n!}\int_{w/2}^{(z+w)/2}\!|p(r)|\,dr\int_0^w\bigg(\frac12\int_0^s\bigg(\int_{r/2}^{(z+r)/2}\!\!|p(u)|\,du\bigg)\,dr\bigg)^n\bigg(\int_{s/2}^{(z+s)/2}\!\!|p(r)|\,dr\bigg)\,ds\en=\\
=\en\frac1{2(n+1)!}\int_{w/2}^{(z+w)/2}\!\!|p(r)|\,dr\,\bigg(\frac12\int_0^w\bigg(\int_{s/2}^{(z+s)/2}\!\!|p(r)|\,dr\bigg)\,ds\bigg)^{n+1}.
\end{multline*}
In this way the inequality
$$
|A(z,w)|\en\leqslant\en\frac12\int_{w/2}^{(z+w)/2}\!\!|p(r)|\,dr\,\,
\exp\bigg(\frac12\int_0^w\bigg(\int_{s/2}^{(z+s)/2}\!\!|p(r)|\,dr\bigg)\,ds\bigg)
$$
is proved which implies the required estimate.

\bigskip

{\bf Definition 1.} \, {\it Operator $\,T\,$ is intertwining {\rm(}transformation operator{\rm)} for the pair of operators $\,\{X,Y\}\,$ iff
\begin{equation}
T\,(X-\lambda I)^{-1}\en=\en(Y-\lambda I)^{-1}T\,,
\label{form7}
\end{equation}
or equivalently $\,\,TD(X)\subset D(Y)\,\,$ and $\,\,TXf=Y\,Tf\,$ for $\,f\in D(X)\,.$
}

\bigskip

The relationship between Volterra type operator appearing in the integral representation of solution $\,s(x,k)\,$ and transformation operator $\,T\,$ for the pair $\,\{L_0,\widehat{L}\}\,$ is established by the following

\bigskip

{\bf Statement 1.} \, {\it Suppose that for a bounded potential $\,p(x)\,$ condition {\rm (\ref{form5})} is satisfied.
Then integral operator
$$
K:\,f(x)\en\longmapsto\en \int_0^x \!K(x,\xi)f(\xi)\,d\xi
$$
is bounded in $\,{\cal H}={\rm L}_2(0,\infty)\,$ and its norm admits the estimate
$$
\|K\|\en\leqslant\en\int_0^{\infty}\!\!r|\,p(r)|\,dr\,\exp\bigg(\int_0^{\infty}\!\!r|\,p(r)|\,dr\bigg)\,.
$$
Moreover operator $\,T=I+K\!:\,{\cal H}\to\widehat{\cal H}_c\,$ is intertwining for the pair $\,L_0\,$ and $\,\widehat{L}\,$ while its image is dense in $\,\widehat{\cal H}_c\,.$
}

\bigskip

{\bf Proof.} By virtue of the inequality
\begin{equation*}
|K(x,\xi)|\en\leqslant\en \rho(x-\xi)\,,\quad \rho(x)\,:=\,C/2\int_{x/2}^{\infty}|p(r)|\,dr\,,
\end{equation*}
resulting from lemma 1, where $\,\displaystyle{C\,=\,\exp\bigg(\int_0^{\infty}\!\!r|p(r)|\,dr\bigg)},\,$ the estimate
\begin{multline*}
\|Kf\|^2\,=\,\int_0^{\infty}\bigg|\!\int_0^x \!K(x,\xi)f(\xi)\,d\xi\,\bigg|^2\!dx\,\leqslant\,\int_0^{\infty}\!\!\bigg(\int_0^x\!\rho(x-\xi)\,d\xi\bigg)\!
\bigg(\int_0^x\!\rho(x-\xi)|f(\xi)|^2d\xi\bigg)dx\\
\leqslant\en \int_0^{\infty}\!\rho(s)\,ds\int_0^{\infty}\!\bigg(\int_{\xi}^{\infty}\!\rho(x-\xi)\,dx\bigg)|f(\xi)|^2d\xi\en=\en \bigg(\int_0^{\infty}\!\rho(s)\,ds\bigg)^2\|f\|^2
\end{multline*}
holds true and therefore one has
\begin{multline*}
\|K\|\en\leqslant\en \int_0^{\infty}\!\rho(s)\,ds\en=\en \int_0^{\infty}\!\!\bigg(\int_x^{\infty}\!\!|p(r)|\,dr\bigg)dx\,\,\exp\bigg(\int_0^{\infty}\!\!r|p(r)|\,dr\bigg)\en=\\
=\en\int_0^{\infty}\!\!r|p(r)|\,dr\,\,\exp\bigg(\int_0^{\infty}\!\!r|p(r)|\,dr\bigg)\,.
\end{multline*}

Let transform $\,\widehat{\Psi}\,$ be defined by the formula
$$
\widehat{\Psi}\,f(k)\en=\en\sqrt{\frac2{\pi}}\int_0^{\infty}\!\frac{s(x,k)}{\hat{e}(0,k)}\,f(x)\,dx\,,
$$
where $\,\hat{e}(x,k)\,$ stands for Jost solution to equation (\ref{form6}).
Eigenfunction expansion theorem for Schr\"odinger operator (see [14]) implies that under condition (\ref{form5}) transform $\,\widehat{\Psi}\!:\,\widehat{\cal H}_c\to{\cal H}\,$ is isometric and an inverse to $\,\widehat{\Psi}\,$ has the form
$$
\widehat{\Psi}^{-1}g(x)\en=\en\sqrt{\frac2{\pi}}\int_0^{\infty}\!\frac{s(x,k)}{\hat{e}(0,-k)}\,g(k)\,dk\,.
$$
Herein $\,\ker\widehat{\Psi}=\widehat{\cal H}_c^{\bot}\,$ and $\,\widehat{\Psi}\,$ implements spectral representation of operator $\,\widehat{L}\,$ in the subspace $\,\widehat{\cal H}_c,\,$ i.e.
\begin{equation}
\big(\widehat{\Psi}\widehat{L}f\big)(k)\en=\en k^2\big(\widehat{\Psi}f\big)(k)\,,\quad f\in D(\widehat{L})\cap\widehat{\cal H}_c\,.
\label{form8}
\end{equation}
Denote by $\,\Phi\,$ a standard Fourier sine-transform
$$
\Phi f(k)\en=\en\sqrt{\frac2{\pi}}\int_0^{\infty}f(x)\,\sin kx\,dx\,,
$$
and let $\,E\,$ be operator of multplication by function $\,\hat{e}(-k):=\hat{e}(0,-k)$ in $\,{\rm L}_2(0,\infty).\,$

\smallskip

One can verify that $\widehat{\Psi}^{-1}E\,\Phi=T.\,$ In fact for arbitrary $\,f\in{\rm C}_0^{\infty}(\mathbb R_+)\,$ one has
\begin{multline*}
\big(\widehat{\Psi}^{-1}E\,\Phi f\big)(x)\en=\en\frac2{\pi}\int_0^{\infty}s(x,k)\bigg(\int_0^{\infty}f(\xi)\sin k\xi\,d\xi\bigg)\,dk\en=\\
=\en\frac2{\pi}\int_0^{\infty}\Phi f(k)\,\sin kx\,dk\en+\en\frac2{\pi}\int_0^{\infty}\!\!\bigg(\int_0^xK(x,\xi)\sin k\xi\,d\xi\bigg)\Phi f(k)\,dk\,.
\end{multline*}
Changing the order of integration in the second summand of the previous line right-hand-side we get
\begin{multline*}
\big(\widehat{\Psi}^{-1}E\,\Phi f\big)(x)\en=\en f(x)\en+\en\frac2{\pi}\int_0^x\!K(x,\xi)\bigg(\int_0^{\infty}\Phi f(k)\,\sin k\xi\,dk\bigg)\,d\xi\en=\\
=\en f(x)\en+\en\int_0^x\!K(x,\xi)f(\xi)\,d\xi\en=\en Tf(x)\,.
\end{multline*}
In this way the required equality is established for $\,f\in{\rm C}_0^{\infty}(\mathbb R_+)\,$ and can be extended to $\,{\rm L}_2(\mathbb R_+)\,$ by continuity.

\smallskip

Now let us show that $\,T=\widehat{\Psi}^{-1}E\,\Phi:\, {\cal H}\to\widehat{\cal H}_c\,$ is a transformation operator for the pair $\,\{L_0,\widehat{L}\}.\,$ Indeed one has
$$
\widehat{\Psi} T (L_0-\lambda I)^{-1}\,\,=\,\,E\Phi(L_0-\lambda I)^{-1}\,\,=\,\, E(k^2-\lambda)^{-1}\Phi\,,
$$
since $\,\Phi\,$ realizes spectral representation of $\,L_0.\,$ On the other hand in virtue of (\ref{form8}) we get
$$
\widehat{\Psi} (\widehat{L}-\lambda I)^{-1}T\,\,=\,\,(k^2-\lambda)^{-1}\widehat{\Psi} T\,\,=\,\,(k^2-\lambda)^{-1}E\Phi
$$
and thus relationship (\ref{form7}) for the pair $\{L_0,\widehat{L}\}$ is established.

\smallskip

Provided that operator $\,\widehat{L}\,$ does not possess a virtual level at the edge of its continuous spectrum (i.e. $\,\hat{e}(0)\ne 0$ ) the corresponding Jost function $\,\hat{e}(k)\,$ is bounded away from zero and consequently operator $\,T=\widehat{\Psi}^{-1}E\,\Phi:\, {\cal H}\to\widehat{\cal H}_c\,$ has a bounded inverse. Anyhow $\,\Phi^{-1}E^{-1}{\rm C}_0^{\infty}(\mathbb R_+)\subset{\rm L}_2(\mathbb R_+)\,$ hence the image $\,T\cal H\,$ contains the set $\,\widehat{\Psi}^{-1}{\rm C}_0^{\infty}(\mathbb R_+)\,$ and thus it is dense in absolutely continuous spectrum subspace $\,\widehat{\cal H}_c\,$ of operator $\,\widehat{L}.$

\bigskip

Below in \S\,6 the proof of density property for the so-called smooth vectors in $\,\widehat{\cal H}_c\,$ (which will be required in the construction of direct and inverse wave operators for the pair $\,\{L,\widehat{L}\})\,$ makes usage of the following preparatory statement about the integral kernel of transformation operator $\,T\,.$

\bigskip

{\bf Lemma 2.} {\it Provided that condition {\rm(\ref{form5})} is satisfied certain constants $\,C_1, C_2, C_3\,$ exist such that given a one-parameter function family $\,\,\psi_a(x,t)\,=\,\psi(x+a,t),\,$ where $\,\,\displaystyle{\psi(x,t)\,=\,x\exp\bigg(\!\!-\frac{x^2}{4(1+it)}\bigg)},\,\,a\in\mathbb R,\,$ the estimate
$$
\int_0^x|K(x,\xi)|\,|\psi_a(\xi,t)|\,d\xi\en\leqslant\en C_1\,x^{1/3}\,+\,\,C_2\,x\exp\bigg(\!\!-\frac{x^{4/3}}{8(1+t^2)}\bigg)\,+\,\,C_3\,x^{4/3}\exp\bigg(\!\!-\frac{x^2}{10(1+t^2)}\bigg)
$$
holds for $\,x\,$ large enough and arbitrary $\,t\geqslant 0.\,$
}

\bigskip

{\bf Proof.} For sufficiently large $\,x\,$ and all $\,t\geqslant 0\,$ in virtue of the estimate $\,\,|K(x,\xi)|\,\leqslant\, \rho(x-\xi)\,\,$ the inequality
$$
\int_0^x\!|\,K(x,\xi)|\,|\,\psi_a(\xi,t)|\,d\xi\en\leqslant\en\bigg(\int_0^{x^{2/3}}\!\!\!+\,\,\int_{x^{2/3}}^{2x/3}\!\!+\,\,\int_{2x/3}^x\bigg)\,\rho(x-\xi)\,|\,\xi+a|\,
\exp\bigg(\!\!-\frac{(\xi+a)^2}{4(1+t^2)}\bigg)\,d\xi\,
$$
is valid. Taking the estimate $\,\,\rho(\eta)\,\leqslant\,C\ln C\big/\eta\,\,$ into account one can separately evaluate each of the summands $\,I_1, I_2, I_3\,$ on the right-hand-side of the latter inequality. As regards the first two of them the constants $\,C_1, C_2\,$ exist such that
\begin{eqnarray*}
I_1&\leqslant&\rho(\,x-x^{2/3})\int_0^{x^{2/3}}(\xi+|a|)\,d\xi\en\leqslant\en C_1\,x^{1/3}\,,\\
I_2&\leqslant&\rho(x/3)\exp\bigg(\!\!-\frac{(x^{2/3}+a)^2}{4(1+t^2)}\bigg)\int_{x^{2/3}}^{2x/3}(\xi+a)\,d\xi\en\leqslant\en C_2\,x\exp\bigg(\!\!-\frac{x^{4/3}}{8(1+t^2)}\bigg)
\end{eqnarray*}
for $\,x\,$ large enough and arbitrary $\,t\geqslant 0.\,$ The third summand $\,I_3\,$ can be represented in the form
$$
I_3\en=\en\bigg(\int_{2x/3}^{x-x^{2/3}}\!\!\!+\,\,\int_{x-x^{2/3}}^{x-x^{1/3}}\!\!+\,\,\int_{x-x^{1/3}}^x\bigg)\,\rho(x-\xi)\,|\,\xi+a|\,
\exp\bigg(\!\!-\frac{(\xi+a)^2}{4(1+t^2)}\bigg)\,d\xi\,.
$$
Each of the integrals here are to be dealt with in the same manner as above in the cases of $\,I_1\,$ or $\,I_2\,$ and thus for sufficiently large $\,x\,$ we get the inequality \begin{multline*}
I_3\en\leqslant\en\rho(x^{2/3})\exp\bigg(\!\!-\frac{(2x/3+a)^2}{4(1+t^2)}\bigg)\int_{2x/3}^{x-x^{2/3}}\!\!(\xi+a)\,d\xi\en+\\
+\en\rho(x^{1/3})\exp\bigg(\!\!-\frac{(x-x^{2/3}+a)^2}{4(1+t^2)}\bigg)\int_{x-x^{2/3}}^{x-x^{1/3}}\!\!(\xi+a)\,d\xi\en+\\
+\en\rho(0)\exp\bigg(\!\!-\frac{(x-x^{1/3}+a)^2}{4(1+t^2)}\bigg)\int_{x-x^{1/3}}^x(\xi+a)\,d\xi\en\leqslant\en C_3\,x^{4/3}\exp\bigg(\!\!-\frac{x^2}{10(1+t^2)}\bigg)\,
\end{multline*}
with a certain constant $\,C_3\,$ which does not depend on $\,t.\,$

\bigskip

\begin{center}
\bf \S\,5. Spectral components of operator $\,L$
\end{center}

\medskip

{\it Characteristic function} for a dissipative Schr\"odinger operator $\,L=L_0+V\,$ with bounded potential $\,V(x)=p(x)+iq(x)\,$ is given (see [18]) by the expression
$$
S(\lambda)\,\,=\,\,I\,+\,2i\,\sqrt{Q}\,\big(L^*-\lambda I\big)^{-1}\!\sqrt{Q}
$$
where $\,Q\,$ denotes an operator of multiplication by $\,q(x).\,$ Analytic in $\,\mathbb C_+\,$ operator valued function $\,S(\lambda)\,$ proves to be a contraction mapping of the space $\,{\cal E}:=\overline{Q\cal H}\,$ and for almost all $\,k\in\mathbb R\,$ possesses boundary limit values $\,S(k):=S(k+i0)\,$ in the sense of strong convergence. Moreover given $\,\lambda\notin\sigma_d(L)\,$ one has
$$
S^{-1}(\lambda)\,\,=\,\,I\,-\,2i\,\sqrt{Q}\,\big(L-\lambda I\big)^{-1}\!\sqrt{Q}\,.
$$

A bounded analytic in $\,\mathbb C_+\,$ function $\,m(\lambda)\not\equiv 0\,$ is called (see [18]) a scalar multiple for $\,S(\lambda)\,$ if analytic in $\,\mathbb C_+\setminus\sigma_d(L)\,$ operator function $\,m(\lambda)S^{-1}(\lambda)\,$ is bounded so that its singularities at points $\,\lambda\in\sigma_d(L)\,$ are removable. A scalar multiple $\,m(\lambda)\,$ being analytic and bounded in $\,\mathbb C_+\,$ admits (see [19]) Nevanlinna-Riesz canonical factorization
$$
m(\lambda)\en=\en m_1(\lambda)\cdot m_2(\lambda)\cdot m_3(\lambda)\,;
$$
here $\,m_1(\lambda)\,$ stands for Blaschke product vanishing at the point set $\,\sigma_d(L),\,$ factor $\,m_2(\lambda)\,$ is a singular inner function of the form
$$
m_2(\lambda)\en=\en e^{ia\lambda}\exp\bigg(i\int_{\mathbb R}\frac{t\lambda-1}{\lambda+t}\,d\mu(t)\bigg)\,,\quad a\geqslant 0\,,
$$
generated by a certain singular measure $\,\mu,\,$ while an outer function $\,m_3(\lambda)\,$ is determined by the boundary values $\,m(k+i0)\,$ according to the formula
$$
m_3(\lambda)\en=\en e^{ib}\exp\bigg(\frac{i}{\pi}\int_{\mathbb R}\ln|\,m(k+i0)|\,\frac{k\lambda+1}{\lambda-k}\,\frac{dk}{1+k^2}\bigg)\,,\en b\in\mathbb R\,.
$$

\medskip

The problem of searching for (and study the structure of) a scalar multiple for the characteristic function $\,S(\lambda)\,$ reduces to investigation of analytic properties of the resolvent $\,R(\lambda)=(L-\lambda I)^{-1}\,$ which for $\,\lambda=k^2,\,{\rm Im}\,k>0,\,$ is an integral operator with the kernel
$$
R(x,\xi,\lambda)\en=\en y\big(\min\{x,\xi\},k\,\big)\,e\big(\max\{x,\xi\},k\,\big)\big/\,e(k)\,,
$$
where $\,y(x,k)\,$ is a solution to equation (\ref{form1}) satisfying the initial conditions $\,y(0,k)=0,$ $\,y'_{x}(0,k)=1.\,$
Denote by $\,A\,$ and $\,B\,$ operators of multiplication by functions $\,a(x)\,$ and $\,b(x)\,$ such that
$$
\left\langle\,a\right\rangle\,:=\,\bigg(\int_0^{\infty}\!x\,|a(x)|^2\,dx\bigg)^{1/2}\!<\en\infty\,,\quad
\left\langle\,b\right\rangle\,:=\,\bigg(\int_0^{\infty}\!x\,|b(x)|^2\,dx\bigg)^{1/2}\!<\en\infty\,.
$$

\medskip

{\bf Lemma 3.} {\it Under the assumption {\rm (\ref{form2})} operator function $\,e(k)AR(k^2)B\,$ is analytic in $\,\mathbb C_+\,$ and for all $\,k\in\mathbb C_+\,$ satisfies the inequality
$$
\|e(k)AR(k^2)B\|\en\leqslant\en K\left\langle\,a\right\rangle\left\langle\,b\right\rangle\,,\quad K=\,2\,\exp\bigg(2\int_0^{\infty}\!x\,|V(x)|\,dx\bigg).
$$
}

\medskip

{\bf Proof.}
For arbitrary $\,\tau:={\rm Im}\,k\geqslant 0,\,k\ne 0,\,$ and all $\,x\in\mathbb R_+\,$ solutions $\,y(x,k)\,$ and $\,e(x,k)\,$ satisfy the estimates
\begin{eqnarray*}
|y(x,k)|&\leqslant&2\,x\exp\bigg(2\int_0^x\xi|V(\xi)|\,d\xi\bigg)\,e^{\tau x}\,,\\
|e(x,k)|&\leqslant&\exp\bigg(\int_x^{\infty}\xi|V(\xi)|\,d\xi\bigg)\,e^{-\tau x}\,,
\end{eqnarray*}
the first of which prove to be an inequality of Gronwall type, while the second one is obtained by the usage of iteration series representation for $\,e(x,k)\,$ (cf. proof of theorem 1). In virtue of the estimates indicated above for the resolvent integral kernel the inequality
$$
|R(x,\xi,\lambda)|\en\leqslant\en\frac{K}{|e(k)|}\,\min\{x,\xi\}
$$
is valid and hence one has
\begin{multline*}
\|e(k)AR(k^2)Bf\|^2\en\leqslant\en K^2\int_0^{\infty}|a(x)|^2\,\bigg|\int_0^\infty\min\{x,\xi\} b(\xi)f(\xi)\,d\xi\,\bigg|^2dx\en\leqslant\\
K^2\int_0^{\infty}|a(x)|^2\,\bigg(\int_0^{\infty} \big(\min\{x,\xi\})^2 |b(\xi)|^2\,d\xi\bigg)\bigg(\int_0^{\infty}|f(\xi)|^2\,d\xi\bigg)\,dx\en
\leqslant\en K^2\left\langle\,a\right\rangle^2\left\langle\,b\right\rangle^2\|f\|^2\,.
\end{multline*}
Resolvent $\,R(k^2)\,$ is known to be meromorphic in $\,\mathbb C_+\,$ so that its poles are located at the zeroes of $\,e(k).\,$ In due turn operator function $\,e(k)AR(k^2)B\,$ is analytic and according to the above estimate bounded in the vicinity of each resolvent pole and therefore has a removable singularity therein.

\bigskip

{\bf Statement 2.} {\it Under the condition {\rm (\ref{form2})}
boundary limit values of $\,S^{-1}(\lambda)\,$ exist {\rm (}in the strong convergence sense{\rm )} almost everywhere on the real line and are bounded provided that $\,e(s)\ne 0\,,\,s\in\mathbb R.\,$ Jost function $\,e(\sqrt{\lambda})\,$ is a scalar multiple for the characteristic function $\,S(\lambda)\,$ of dissipative operator $\,L.\,$ If $\,L\,$ has no spectral singularities then the singular factor in the canonical factorization of $\,e(\sqrt{\lambda})\,$ is trivial.
}

\medskip

{\bf Proof} employs the estimate
$$
\big\|\,e(k)\sqrt{Q}\,R(k^2)\sqrt{Q}\,\big\|\,\,\leqslant\,\,K\langle\sqrt{q}\,\rangle^2
$$
established in lemma 3 to verify that analytic in $\,\mathbb C_+\,$ operator function $\,e(k)\sqrt{Q}\,R(k^2)\sqrt{Q}\,$ possesses boundary limit values almost everywhere on $\,\mathbb R.\,$ Since $\,e(k)\,$ is continuous up to the real axis and $\,e(s)\ne 0\,$ for $\,s\in\mathbb R\,$ it follows that (strong convergence) limit values of operator function $\,S^{-1}(\lambda)\,$ exist and are bounded almost everywhere on the real line as well. Moreover operator function $\,e(\sqrt{\lambda})\,S^{-1}(\lambda)\,$ is analytic and bounded in $\,\mathbb C_+\,$
including the points $\,\lambda\in\sigma_d(L),\,$ so that $\,m(\lambda)=e(\sqrt{\lambda})\,$ proves to be a scalar multiple for the characteristic function $\,S(\lambda)\,$ of operator $\,L.\,$ In the absence of spectral singularities discrete spectrum $\,\sigma_d(L)\,$ is finite and hence factor $\,m_1(\lambda)\,$ is a finite Blaschke product. Within the above context we shall make use of the following assertion (see [19]) from analytic functions theory: together with $\,m(\lambda)\,$ its outer factor $\,m_3(\lambda)\,$ is continuous up to the real line so that $\,|m(k+i0)|=|m_3(k+i0)|,\,$ while the singular inner function $\,m_2(\lambda)\,$ can be extended by continuity to only  points of the real axis which does not belong to the support of measure $\,\mu.\,$ Provided that operator $\,L\,$ has no spectral singularities $\,m_2(\lambda)\,$ proves to admit an extension by continuity on entire $\,\mathbb R,\,$ hence $\,{\rm supp}\,\mu=\varnothing\,$ and therefore $\,m_2(\lambda)\equiv 1.$

\bigskip

Further by means of Riesz integral define a projection
$$
P\,\,=\,\,-\,\frac1{2\pi i}\oint_{\Gamma}\,(L-\lambda I)^{-1}\,d\lambda
$$
onto the linear span $\,{\cal H}_d:=P{\cal H}\,$ of root (eigen and associated) vectors of operator $\,L\,$ parallel to the subspace  $\,{\cal H}_c:=(I-P){\cal H};\,$ here $\,\Gamma\subset\mathbb C_+\,$ is appropriately oriented closed contour separating spectral components $\,\sigma_c(L)\,$ and $\,\sigma_d(L).\,$ Operator $\,L\,$ can be decomposed with respect to the direct sum representation
$\,{\cal H}\,=\,{\cal H}_c\,\dot{+}\,{\cal H}_d\,$ in the following (see e.g. [20]) sense:
$$
PD(L)\subset D(L)\,,\quad L{\cal H}_c\subset{\cal H}_c\,,\quad L{\cal H}_d\subset{\cal H}_d\,,
$$
while $\,\sigma(L|{\cal H}_c)=\sigma_c(L)\,$ and $\,\sigma(L|{\cal H}_d)=\sigma_d(L).$

\bigskip

{\bf Definition 2.} {\it Absolutely continuous {\rm (}outer{\rm )} subspace} ${\cal N}_e\,$ for completely nonselfadjoint dissipative operator $\,L\,$ is the closure in $\,\cal H\,$ of the set of elements $\,\psi\in\cal H\,$ for which vector function $\,\sqrt{Q}R(\lambda)\psi\,$ is analytic in $\,\mathbb C_+\,$ and
$$
\sup_{\varepsilon>0}\,\int_{\mathbb R}\big\|\sqrt{Q}R(k+i\varepsilon)\psi\,\big\|^2\,dk\,\,<\,\,\infty\,.
$$
{\it Singular {\rm (}inner{\rm )} subspace} $\,{\cal N}_i\,$ of operator $\,L\,$ consists of vectors $\,\psi\in{\cal H}\,$ such that weak limit values $\,R(k\pm i0)\psi\,$ coincide for almost all $\,k\in\mathbb R.$

\bigskip

A direct sum decomposition
$$
{\cal H}\,=\,{\cal N}_e\,\dot{+}\,{\cal N}_i\,,\,\en {\cal N}_e\subset{\cal H}_c\,,\,\en
{\cal H}_d\subset{\cal N}_i\,,
$$
takes place and moreover projection on $\,{\cal N}_i\,$ parallel to $\,{\cal N}_e\,$ is bounded provided that operator $\,L\,$ has no spectral singularities (see [21]). Below we shall need certain information about analytic in a disc contractive operator valued functions (see [18]). By $\,{\rm H}_2({\mathfrak A})\,$ we denote Hardy space of analytic in the disc $\,{\mathbb D}=\{|z|<1\}\,$ vector valued in a Hilbert space $\,\mathfrak A\,$ functions which are identified with equivalence classes of their boundary values on the circle $\,|z|=1,\,$ and in this way treated as a subspace in $\,{\rm L}_2({\mathfrak A}).\,$ An analytic in the disc $\,\mathbb D\,$ contractive operator function $\,\Theta(z)\,$ acting from Hilbert space $\,\mathfrak A\,$ into Hilbert space $\,\mathfrak A'\,$ is called {\it outer}, if $\,{\rm clos}_{{\rm L}^2({\mathfrak A'})}\Theta{\rm H}_2({\mathfrak A})={\rm H}_2({\mathfrak A'}).\,$ An operator function  $\,\Theta(z):\mathfrak A\to\mathfrak A'\,$ is called {\it inner} if the mapping $\,\Theta\,$ from $\,{\rm H}_2({\mathfrak A})\,$ to $\,{\rm H}_2({\mathfrak A'})\,$ proves to be isometric. Any analytic in $\,\mathbb D\,$ contractive operator function $\,\Theta(z)\,$ can be represented in the form $\,\Theta(z)=\widehat{\Theta}(z)\oplus{\rm const}\,$ where $\,\widehat{\Theta}(z)\,$ is purely contractive (i.e. $\|\widehat{\Theta}(0)\|<1)$ while the constant summand is unitary.

\bigskip

{\bf Assertion 1.} {\it Given completely nonselfadjoint dissipative operator $\,L,\,$ for which $\,\sigma_c(L)\subset\mathbb R\,$ and $\,\#\,\sigma_d(L)<\infty,\,$ one has $\,{\cal N}_i\cap{\cal H}_c=\{0\},\,$ provided that the singular factor in the canonical factorization of scalar multiple for the corresponding characteristic function $\,S(\lambda)\,$ is trivial.
}

\bigskip

{\bf Proof.} Denote by $\,\,C:=(L-iI)(L+iI)^{-1}\,$ Cayley transform of a dissipative completely nonselfadjoint operator $\,L\,$ and consider its characteristic function
$$
\Theta(z)\,\,=\,\,-\,C\,+\,z\,(I-CC^*)^{1/2}(I-z\,C^*)^{-1}\,(I-C^*C)^{1/2},
$$
being a contraction mapping from $\,{\mathfrak A}={\rm clos}_{\cal H}(I-C^*C)^{1/2}{\cal H}\,$ to $\,{\mathfrak A'}={\rm clos}_{\cal H}(I-CC^*)^{1/2}{\cal H},\,$ which is analytic in $\,z\in\mathbb D\,$ and coincides with $\,S\big(i(1+z)/(1-z)\big).\,$ Thereby operator function $\,\Theta(z)\,$ possesses a scalar multiple such that inner factor in its Nevanlinna-Riesz canonical factorization reduces to a finite Blaschke product.

Operator function $\,\Theta(z):\,\mathfrak A\to\mathfrak A'\,$ admits (see [18]) a factorization of the form $\,\,\Theta(z)=\Theta_1(z)\Theta_2(z),\,\,$
where $\,\Theta_1^*(\bar z)\,$ is an outer operator valued function, while $\,\Theta_2^*(\bar z)\,$ is an inner one. By Nagy and Foias theorem on inner-outer factorization of operator functions possessing a scalar multiple, Blaschke product
$$
\beta(z)\,\,=\,\,\prod_n\,\frac{z-z_n}{1-\bar z_nz}\,,\quad z_n=\frac{\lambda_n-i}{\lambda_n+i}\,,\quad\lambda_n\in\sigma_d(L),
$$
proves to be a scalar multiple for the factor  $\,\Theta_2(z).\,$
Due to the theorem on triangulation of a contraction generated by a factorization of its characteristic function (see [18]) the subspace  $\,{\cal N}_i\,$ is invariant with respect to $\,C\,$ and characteristic function of the restriction $\,C|{\cal N}_i\,$ coincides with the pure part $\,\widehat{\Theta}_2(z).$

The subspace $\,{\cal N}_s:={\cal N}_i\,\cap\,{\cal H}_c=(I-P){\cal N}_i\,$ is also invariant under the action of $\,C\,$ and $\,{\cal N}_i={\cal H}_d\,\dot{+}\,{\cal N}_s.\,$ The indicated decomposition generates a factorization $\,\,\widehat{\Theta}_2(z)=\Theta_3(z)\Theta_4(z)\,\,$ such that characteristic function of the restriction $\,C|{\cal N}_s\,$ coincides with $\,\widehat{\Theta}_4(z).\,$ Since $\,\#\,\sigma_d(L)<\infty\,$ the subspace $\,{\cal N}_s\,$ in $\,{\cal N}_i\,$ has a finite codimension and hence $\,\Theta_3(z)\,$ proves to be a contracting operator function in a finite dimensional space certainly admitting a scalar multiple. By this reason factor $\,\Theta_4(z)\,$ possesses the same scalar multiple as $\,\widehat{\Theta}_2(z),\,$ namely $\,\beta(z),\,$ and so operator function $\,\widehat{\Theta}_4(z)\,$ with such a property is specified to be inner.

Taking into account that $\,\sigma(L_c)\subset\mathbb R\,$ one has $\,\sigma(C|{\cal H}_c)\cap\mathbb D=\varnothing\,$ and therefore characteristic function of the restriction $\,C|{\cal N}_s\,$ is invertible in $\,\mathbb D.$ Thus operator function $\,\widehat{\Theta}_4(z)\,$ is invertible for any $\,z\in\mathbb D\,$ and possesses a scalar multiple $\,\beta(z).$ It follows that  $\,\widehat{\Theta}_4^{-1}(z)\,$ is analytic and uniformly bounded in $\,\mathbb D\,$ and hence $\,\widehat{\Theta}_4(z)\,$ proves to be outer. An operator function inner and outer at the same time is known to be a unitary constant and as applied to $\,\widehat{\Theta}_4(z)\,$ such a specification implies that $\,{\cal N}_s=\{0\}.$

\bigskip

Combining assertion 1 with statement 2 one gets

\bigskip

{\bf Statement 3.} {\it If a completely nonselfadjoint dissipative operator $\,L=L_0+V\,$ with a potential $\,V(x)\,$ satisfying condition  {\rm(\ref{form2})} has no spectral singularities then $\,\,{\cal N}_e={\cal H}_c\,\,$ and $\,\,{\cal H}_d={\cal N}_i\,.$
}

\bigskip

In conclusion of this section we give a short summary of the required facts (see [22]) about functional model of a completely nonselfadjoint dissipative operator having $\,L=\widehat{L}+iQ\,$ in mind to be dealt with. By $\,\mathfrak L\,$ denote a Hilbert factor-space of two-component vector functions taking values in $\,\cal E\oplus\cal E\,$ which are supposed to be square-integrable on $\,\mathbb R\,$ with the operator weight
$$
\begin{pmatrix}
I&S^*(k)\\
S(k)&I
\end{pmatrix}\,.
$$
In $\,\mathfrak L\,$ define a subspace
$$
\mathfrak H\,\,=\,\,\bigg\{\begin{pmatrix}f\\g\end{pmatrix}\in\mathfrak L:\,
f+S^*g\in{\rm H}_2^-({\cal E}),\, Sf+g\in{\rm H}_2^+({\cal E})\bigg\}\,;
$$
here $\,{\rm H}_2^{\pm}(\cal E)\,$ denote Hardy classes of boundary values attained on the real line by analytic vector functions  $\,h:\mathbb C_{\pm}\to\cal\, E\,$ for which
$$
\sup_{\varepsilon>0}\,\int_{\mathbb R}\|h(k\pm i\varepsilon)\|^2\,dk\,<\,\infty\,.
$$
Completely nonselfadjoint dissipative operator $\,L\,$ is equivalent to the generator of contraction semigroup
$$
Z(t)\,\,=\,\,{\cal P}\,{\cal U}(t)|{\mathfrak H}\,,\quad t\geqslant 0\,,
$$
where unitary group $\,{\cal U}(t)\,$ acts in $\mathfrak L\,$ by the formula
$$
{\cal U}(t)\begin{pmatrix}f\\g\end{pmatrix}\!(k)\,\,=\en
\exp(ikt)\begin{pmatrix}f(k)\\g(k)\end{pmatrix}\,,
$$
while $\,\cal P\,$ is an orthoprojection from $\mathfrak L\,$ onto $\mathfrak H.\,$ Generator of the semigroup $\,Z(t)\,$ is called a {\it functional model} of operator $\,L.\,$ Isometry implementing the equivalence indicated above will be denoted by $\,J:{\cal H}\to{\mathfrak H}.\,$

\bigskip

\begin{center}
\bf \S\,6. Direct and inverse wave operators
\end{center}

\medskip

Let us put into context two unitary groups $\,U_0(t)=\exp(itL_0)\,$ and $\,\widehat{U}(t)=\exp(it\widehat{L})\,$ and a one-parameter operator family
$$
U(t)\,=\,\exp(itL)\,=\,{\rm s}\!\!-\!\!\!
\lim_{n\to\infty}\big(I-itL/n\big)^{-n}.
$$

The goal of this section is to prove existence and completeness of direct and inverse wave operators for the pair $\,\{L,\widehat{L}\}.\,$ Beforehand we shall establish the density of the so-called {\it smooth vectors} in $\,\widehat{\cal H}_c\,$ which will be used below in the construction of wave operators.

\bigskip

{\bf Lemma 4.} {\it Provided that integral {\rm (\ref{form2})} converges subspace $\,\widehat{\cal H}_c\,$ contains a dense subset of elements  $\,\varphi\,$ for which the following property
\begin{equation}
\sqrt{Q}\,\big(\widehat{L}-(k\pm i0) I\big)^{-1}\!\varphi\,\in\,{\rm H}_2^{\pm}(\cal E)
\label{form9}
\end{equation}
is valid.
}

\smallskip

{\bf Proof.} In virtue of the known equalities
$$
\int_{\mathbb R\pm i\varepsilon}\big\|\sqrt{Q}\,(\widehat{L}-\lambda I)^{-1}\varphi\big\|^2d\lambda\,\,=\,\,
2\pi\int_0^{\infty}e^{-2\varepsilon t}\big\|\sqrt{Q}\,\widehat{U}(\mp t)\,\varphi\big\|^2dt\,,
$$
it clearly suffices to establish the existence of the set dense in $\,\widehat{\cal H}_c\,$ elements $\,\varphi\,$ of which satisfy the requirement $\,\big\|\sqrt{Q}\,\widehat{U}(\pm t)\,\varphi\,\big\|\in{\rm L}_2(\mathbb R_+).\,$
Note that $\,{\rm L}_2(0,\infty)\,$ contains a dense subset formed by linear combinations of functions
$$
\phi_a(x)\en=\en \Phi^{-1}\big(\,k e^{-k^2}\cos ak\big)(x)\,,\quad a\in\mathbb R\,,
$$
therefore linear combinations of functions $\,\varphi_a(x)\,=\,T\phi_a(x)\,$ are dense in $\,T\cal H\,$ and hence in $\,\widehat{\cal H}_c$ (cf. statement 1).

\medskip

Taking into account that $\,T\,$ is a transformation operator for the pair $\,L_0\,$ and $\,\widehat{L}\,$ while $\,\Phi\,$ implements spectral representation of operator $\,L_0\,$ we get
$$
\widehat{U}(-t)\varphi_a(x)\en=\en \widehat{U}(-t)T\phi_a(x)\en=\en TU_0(-t)\phi_a(x)
$$
and besides $\,\,\big(\Phi\,U_0(-t)\phi_a\big)(k)\,=\, k\,e^{-k^2(1+it)}\cos ak.\,$
Thus
\begin{multline*}
\widehat{U}(-t)\varphi_a(x)\en=\en T\,\Phi^{-1}\big(k\,e^{-k^2(1+it)}\cos ak\big)(x)\en=\\
=\en\,\frac{(1+it)^{-3/2}}{4\sqrt{\pi}}\,\bigg(\,\psi_a(x,t)\,+\,\psi_{-a}(x,t)\,+\int_0^x K(x,\xi)\big(\psi_a(\xi,t)+\psi_{-a}(\xi,t)\big)\,d\xi\bigg)
\end{multline*}
where $\,\psi_{\pm a}(x,t)\,=\,\psi(x\pm a,t)\,$ and therefore one has
\begin{multline*}
\big\|\sqrt{Q}\,\widehat{U}(\pm t)\,\varphi_a\big\|^2\,\leqslant\,\, \frac1{16\pi}\,(1+t^2)^{-3/2}\bigg[\int_0^{\infty}
q(x)|\psi_a(x,t)|^2dx\,+\,\int_0^{\infty}q(x)|\psi_{-a}(x,t)|^2dx\,\,\,+\\
\int_0^{\infty}q(x)\bigg(\int_0^x|K(x,\xi)|\,|\psi_a(\xi,t)|\,d\xi\bigg)^2dx\,+\,
\int_0^{\infty}q(x)\bigg(\int_0^x|K(x,\xi)|\,|\psi_{-a}(\xi,t)|\,d\xi\bigg)^2dx\,\bigg]\,.
\end{multline*}

\medskip

In virtue of lemma 2 (which can be applied due to the finiteness of the first momentum condition imposed on the potential) for the integrals
\begin{gather*}
\int^{\infty}(1+t^2)^{-3/2}\bigg(\int_0^{\infty}q(x)|\psi_a(x,t)|^2dx\bigg)\,dt\,,\quad a\in\mathbb R\,,\\
\int^{\infty}(1+t^2)^{-3/2}dt\int_0^{\infty}q(x)\bigg(\int_0^x|K(x,\xi)|\,|\psi_a(\xi,t)|\,d\xi\bigg)^2dx
\end{gather*}
the desired convergence is guaranteed by the following requirements
\begin{eqnarray*}
\int^{\infty}(1+t^2)^{-3/2}dt\int_0^{\infty}q(x)\,x^{2/3}\,dx&<&\infty\,,\\
\int^{\infty}(1+t^2)^{-3/2}dt\int_0^{\infty}q(x)\,x^2\exp\bigg(\!-\frac{x^{4/3}}{4(1+t^2)}\bigg)\,dx&<&\infty\,,\\
\int^{\infty}(1+t^2)^{-3/2}dt\int_0^{\infty}q(x)\,x^{8/3}\exp\bigg(\!-\frac{x^2}{5(1+t^2)}\bigg)\,dx&<&\infty\,.
\end{eqnarray*}
The first of the listed integrals under the hypothesis in question clearly is convergent. Consider now the second one
\begin{multline*}
\int_1^{\infty}(1+t^2)^{-3/2}dt\int_0^{\infty}q(x)\,x^2\exp\bigg(\!-\frac{x^{4/3}}{4(1+t^2)}\bigg)\,dx\en\leqslant\\
\leqslant\en\int_1^{\infty}t^{-3}\,dt\int_0^{t^2}q(x)\,x^2\,dx\en+\en\int_1^{\infty}(1+t^2)^{-3/2}dt\int_{t^2}^{\infty}q(x)\,x^2\exp\bigg(\!-\frac{x^{4/3}}{4(1+t^2)}\bigg)\,dx\,.
\end{multline*}
Changing the order of integration in the first summand of the latter inequality right-hand-side one has
$$
\int_1^{\infty}t^{-3}\,dt\int_0^{t^2}q(x)\,x^2\,dx\en=\en\frac12\int_0^1q(x)\,x^2\,dx\,\,+\,\,\frac12\int_1^{\infty}q(x)\,x\,dx\en<\en\infty\,.
$$
To evaluate the second summand taking into account that $\,\,\displaystyle{x\,\exp\bigg(\!-\frac{x^{4/3}}{4(1+t^2)}\bigg)\,\leqslant\,1}\,\,$ for $\,x\geqslant t^2\,$ and sufficiently large $\,t\,$ we get
$$
\int_1^{\infty}(1+t^2)^{-3/2}dt\int_{t^2}^{\infty}\!q(x)\,x^2\exp\bigg(\!-\frac{x^{4/3}}{4(1+t^2)}\bigg)\,dx\,\leqslant\,\int_1^{\infty}(1+t^2)^{-3/2}dt\int_0^{\infty}\!q(x)\,x\,dx\,<\,\infty\,.
$$
In the same manner the third integral can be estimated:
\begin{multline*}
\int_1^{\infty}(1+t^2)^{-3/2}dt\int_0^{\infty}q(x)\,x^{8/3}\exp\bigg(\!-\frac{x^2}{5(1+t^2)}\bigg)\,dx\en\leqslant\\
\leqslant\en\int_1^{\infty}t^{-3}\,dt\int_0^{t^{6/5}}q(x)\,x^{8/3}\,dx\en+\en\int_1^{\infty}(1+t^2)^{-3/2}dt\int_{t^{6/5}}^{\infty}q(x)\,x^{8/3}\exp\bigg(\!-\frac{x^2}{5(1+t^2)}\bigg)\,dx\,,
\end{multline*}
where
$$
\int_1^{\infty}t^{-3}\,dt\int_0^{t^{6/5}}q(x)\,x^{8/3}\,dx\en=\en\frac12\int_0^1q(x)\,x^{8/3}\,dx\,\,+\,\,\frac12\int_1^{\infty}q(x)\,x\,dx\,<\,\infty
$$
and moreover
$$
\int_1^{\infty}\!\!(1+t^2)^{-3/2}dt\int_{t^{6/5}}^{\infty}\!q(x)\,x^{8/3}\exp\bigg(\!-\frac{x^2}{5(1+t^2)}\bigg)\,dx\,\leqslant\,\int_1^{\infty}\!\!(1+t^2)^{-3/2}dt\int_0^{\infty}\!q(x)\,x\,dx\,<\,\infty\,,
$$
since $\,\,\displaystyle{x^{5/3}\,\exp\bigg(\!-\frac{x^2}{5(1+t^2)}\bigg)\,\leqslant\,1}\,\,$ for $\,x\geqslant t^{6/5}\,$ provided that  $\,t\,$ is large enough.

\bigskip

The above proof makes explicit usage of the finiteness of the first momentum condition for the imaginary part of potential $\,V(x)\,$ while hypothesis (\ref{form5}) draws lemma 2 into the present context.
Note that restriction of the type (\ref{form9}) appears in scattering theory for perturbations relatively smooth in the sense of Kato when wave operators are constructed for selfadjoint operators and even certain nonselfadjoint ones as well (see [6] and [23]).

\smallskip

Further by means of the passage to functional model representation (cf. \S\,5) we establish

\bigskip

{\bf Statement 4.} {\it Suppose that condition {\rm(\ref{form2})} is satisfied and a completely nonselfadjoint operator $\,L=L_0+V\,$ has no spectral singularities. Then an operator
$$
\Omega[\widehat{L},L]\en=\en{\rm s}\!-\!\!\lim_{t\to\infty}\widehat{U}(t)\,U(-t)
$$
is well-defined and bounded on the subspace $\,{\cal N}_e,\,$ it proves to be left inverse for the direct wave operator
$$
\Omega[L,\widehat{L}]\en=\en{\rm s}\!-\!\!\lim_{t\to\infty}U(t)\,\widehat{U}(-t)
$$
with the domain $\,\widehat{\cal H}_c\,$ and at the same time right inverse to $\,\Omega[L,\widehat{L}]\,$ on the absolutely continuous subspace of operator $\,L.$
}

\bigskip

{\bf Proof\,} makes use of the appropriately adapted certain fragments of the wave operators construction elaborated and somewhat adjusted in [24] and [25] respectively.

\smallskip

The direct wave operator for the pair $\,\{L,\widehat{L}\}\,$ when written in the functional model representation of operator $\,L\,$ is given on the set of smooth vectors by the formula
$$
J\,\Omega[L,\widehat{L}]\,J^{-1}\,{\cal P}\begin{pmatrix}f\\g\end{pmatrix}\,\,=\en
{\cal P}\begin{pmatrix}f\\-Sf\end{pmatrix}.
$$
Indeed on its proper domain
$$
{\mathfrak D}\en=\en{\cal P}\,\bigg\{\begin{pmatrix}f\\g\end{pmatrix}\in{\mathfrak L}:\,
f+S^*g+Sf+g=0\,\bigg\}
$$
operator
$$
W\,:\,\,{\cal P}\begin{pmatrix}f\\g\end{pmatrix}\,\,\mapsto\en
{\cal P}\begin{pmatrix}f\\-Sf\end{pmatrix}
$$
can be realized in the form $\,W=\,{\rm s}$-$\lim\limits_{t\to\infty}J\,U(t)\,\widehat{U}(-t)\,J^{-1}\,$ and moreover according to lemma 4 the set
$$
J^{-1}{\mathfrak D}\,\,=\,\,\big\{\varphi\in{\cal H}:\,\sqrt{Q}\,\big(\widehat{L}-(k\pm i0) I\big)^{-1}\varphi\,\in\,{\rm H}_2^{\pm}(\cal E)\,\big\}
$$
is dense in $\,\widehat{\cal H}_c\,$ provided that integral (\ref{form2}) converges. Thus a bounded operator $\,W\,$ can be extended by continuity from $\,\mathfrak{D}\,$ to $\,J\widehat{\cal H}_c\,$ and this continuation coincides with $\,J\,\Omega[L,\widehat{L}]\,J^{-1}.\,$ As a consequence one gets the inclusion $\,\Omega[L,\widehat{L}]\widehat{\cal H}_c\,\subset\,{\cal N}_e\,$ since
$$
J{\cal N}_e\,\,=\,\,\overline{\bigg\{{\cal P}\begin{pmatrix}f\\-Sf\end{pmatrix},\,
\int_{\mathbb R}\big\|(I-S^*S)^{1/2}f(k)\big\|^2\,dk<\infty\bigg\}}\,.
$$

\medskip

By virtue of lemma 3 and provided that integral (\ref{form2}) converges the estimate
$$
\big\|\,\hat{e}(k)\sqrt{Q}\,(\widehat{L}-k^2I)^{-1}\sqrt{Q}\,\big\|\,\,\leqslant\,\,K\langle\sqrt{q}\,\rangle^2
$$
is valid; due to this fact operator function $\,I+S(\lambda)\,$ possesses a scalar multiple $\,\hat{e}(\sqrt{\lambda})\,$ so that analytic in  $\,\mathbb C_+\,$ operator function
$$
\big(I+S(\lambda)\big)^{-1}\,=\,\,\,\frac12\,\Big(I\,-\,i\,\sqrt{Q}\,\big(\widehat{L}-\lambda I\big)^{-1}\!\sqrt{Q}\,\Big)
$$
proves to be bounded when $\,\lambda\,$ is separated away from $\,\sigma_d(\widehat{L})\cup\{0\}\,$ and hence it has strong convergence boundary limit values almost everywhere on the real axis. This enables one to define in a model representation an operator
$$
\widetilde{W}\,:\,\,{\cal P}\begin{pmatrix}f\\-Sf\end{pmatrix}\,\,\mapsto\en{\cal P}\begin{pmatrix}f\\-(I+S^*)^{-1}(I+S)\,f\end{pmatrix}
$$
with the domain dense in $\,J{\cal N}_e\,$ on which it can be expressed in the form $\,\widetilde{W}=$ $\,{\rm s}$-$\lim\limits_{t\to\infty}J\widehat{U}(t)\,U(-t)J^{-1}.$

\medskip

Since the operator $\,L\,$ has no spectral singularities the boundary limit values $\,S(k+i0)^{-1}\,$ according to statement 2 exist almost everywhere on the real line and are bounded. Owing to this fact the group $\,U(t)\,$ is uniformly bounded on the absolutely continuous subspace $\,{\cal N}_e.\,$ Thus operator $\,J^{-1}\widetilde{W}\,J\,$ is bounded and admits a continuation to $\,{\cal N}_e\,$ that enables us to define correctly the inverse wave operator
$$
\Omega[\widehat{L},L]\,\varphi\,\,=\,\,\lim_{t\to\infty}\widehat{U}(t)\,U(-t)\,\varphi\,,\quad\varphi\in{\cal N}_e\,.
$$
The definitions of $\,W\,$ and $\,\widetilde{W}\,$ indicate that operator $\,\Omega[\widehat{L},L]\,$ constructed above proves to be left inverse to the direct wave operator $\,\Omega[L,\widehat{L}]\,$ on $\,\widehat{\cal H}_c\,$ and at the same time its right inverse on the subspace $\,{\cal N}_e\,$ and so $\,\,\Omega[L,\widehat{L}]\,\widehat{\cal H}_c\,=\,{\cal N}_e\,.$ As regards existence and boundedness of direct wave operator note that it makes no difference whether operator $\,L\,$ has any spectral singularities or not.

\bigskip
\bigskip

\begin{center}
\bf \S\,7. Absolutely continuous spectrum of operator $\,L\,$
\end{center}

\medskip

Existence of wave operators $\,\Omega[L,L_0]\,$ and $\,\Omega[L_0,L]\,$ and completeness property established below can be regarded as the absence of singular component in the continuous spectrum of $\,L\,$ in the sense that $\,L_c\,$ proves to be similar to operator $\,L_0\,$ having purely absolutely continuous spectrum.

\bigskip

{\bf Proof of theorem 2.} According to statement 4 the direct wave operator $\,\,\Omega[L,\widehat{L}]\,\,=$ $\,\,{\rm s}$-$\lim\limits_{t\to\infty}U(t)\,\widehat{U}(-t)\,\,$ exists and is bounded on $\,\widehat{\cal H}_c\,$ provided that condition    (\ref{form2}) is satisfied. If additionally $\,L\,$ has no spectral singularities then $\,\,{\cal N}_e={\cal H}_c\,\,$ by virtue of statement 3 and besides operator $\,\,\Omega[\widehat{L},L]\,\,=$ $\,\,{\rm s}$-$\lim\limits_{t\to\infty}\widehat{U}(t)\,U(-t)\,\,$ is well-defined and bounded on $\,{\cal N}_e.$

Wave operator $\,\Omega[L,\widehat{L}]\,$ intertwines $\,\widehat {L}\,$ and $\,L\,:$
$$
L\,\Omega[L,\widehat{L}]\en=\en\Omega[L,\widehat{L}]\,\widehat{L}\,,
$$
so that $\,\Omega[L,\widehat{L}]\big(D(\widehat{L})\cap\widehat{\cal H}_c\big)=D(L)\cap{\cal H}_c\,$ and therefore
$$
L_c\en=\en\Omega[L,\widehat{L}]\,\,\widehat{L}_c\,\,\Omega[\widehat{L},L]\,.
$$

Under the assumption $\,p(x)\in{\rm L}_1(\mathbb R_+)\,$ a complete wave operator $\,\Omega[\widehat{L},L_0]\,$ exists (see e.g. [20]) and hence  $\,\widehat{L}_c\,$ and $\,L_0\,$ are unitarily equivalent:
$$
\widehat{L}_c\en=\en\Omega[\widehat{L},L_0]\,\,L_0\,\,\Omega[L_0,\widehat{L}]\,.
$$

Making use of the chain rule one has
\begin{eqnarray*}
\Omega[L,\widehat{L}]\,\,\Omega[\widehat{L},L_0]&=&\Omega\,:\,\,\,{\cal H}\,\to\,\widehat{\cal H}_c\,\to\,{\cal H}_c\\
\Omega[L_0,\widehat{L}]\,\,\Omega[\widehat{L},L]&=&\widetilde{\Omega}\,:\,{\cal H}_c\,\to\,\widehat{\cal H}_c\,\to\,{\cal H}
\end{eqnarray*}
and thus
$$
L_c\en=\en\Omega\,L_0\,\widetilde\Omega\,.
$$
Operator $\,\widetilde\Omega\,$ being left inverse to $\,\Omega\,$ at the same time proves to be its right inverse on the subspace $\,{\cal H}_c.\,$ In this way we have established that wave operators $\,\Omega\,$ and $\,\widetilde\Omega\,$ implement spectral representation of completely nonselfadjoint dissipative operator $\,L\,$ without spectral singularities.

\medskip

A simple sufficient condition for a dissipative operator $\,L\,$ to be completely nonselfadjoint is given by the following

\bigskip

{\bf Assertion 2.} {\it If an interval $\,\Delta\subset\mathbb R_+\,$ exists such that $\,q(x)>0\,$ for $\,x\in\Delta,\,$ then operator  $\,L=L_0+V\,$ is completely nonselfadjoint.
}

\bigskip

{\bf Proof.} A dissipative operator $\,L=\widehat{L}+iQ\,$ is known to be completely nonselfadjoint on the reducing (invariant) subspace
$$
{\cal L}\en=\en\overline{{\rm Lin}\big\{\big(\widehat{L}-\lambda I\big)^{-1}{\cal E}\,,\,\lambda\in\mathbb C\setminus\mathbb R\big\}}\,.
$$
Since $\,q(x)>0\,$ on the interval $\,\Delta\,$ then for an arbitrary function $\,f\in{\rm C}_0^{\infty}(\Delta)\subset\cal E\,$ given  $\,t>0\,$ and $\,\varepsilon>0\,$ one has
$$
-\,\frac1{2\pi i}\int_{\mathbb R-i\varepsilon}e^{i\lambda t}\big(\widehat{L}-\lambda I\big)^{-1}\!f\,d\lambda\,\,=\,\,\widehat{U}(t)f\,\in\,\cal L
$$
and hence $\,(\widehat{U}(t)f,g)=0\,$ for $\,t>0\,$ provided that $\,g\bot\,\cal L.\,$ Denote by $\,\varphi\,$ and $\psi\,$ respectively the images of $\,f\,$ and $\,g\,$ under the action of $\widehat{\Psi}$-transform. Then $\,\varphi\,\overline{\psi}\in{\rm L}_1(\mathbb R_+)\,$ and besides function
$$
\sqrt{\frac{\pi}2}\,\hat{e}(k)\varphi(k)\,=\,\int_{\Delta}\!f(x)s(x,k)\,dx\,=\,\int_{\Delta}\!f(x)\sin kx\,dx\,+\,\int_{\Delta}\!f(x)\bigg(\int_0^xK(x,\xi)\sin k\xi\,d\xi\bigg)dx
$$
is analytic. Since $\,\widehat{\Psi}\,$ realizes spectral representation of operator $\,\widehat{L}_c\,$ we have
$$
(\widehat{U}(t)f,g)\en=\en\int_0^{\infty}\!e^{ik^2t}\varphi(k)\overline{\psi}(k)\,dk\,\,+\,\,\vartheta(t)\,,
$$
where $\,\vartheta(t)\,$ is a quasi-periodic function being a finite sum of terms $\,e^{i\lambda t}(f,h)\,$ such that $\,h\in\widehat{\cal H}_c^{\bot}\,$ and $\,\lambda\in\sigma_d(\widehat{L}).\,$ Thus condition $\,(\widehat{U}(t)f,g)=0\,$ may be satisfied for any $\,f\in{\rm C}_0^{\infty}(\Delta)\,$ and every $\,t>0\,$ iff $\,g\in\widehat{\cal H}_c\,$ so that $\,\vartheta(t)\equiv 0\,$ and
$$
\int_0^{\infty}\!e^{ik^2t}\varphi(k)\overline{\psi}(k)\,dk\,\,=\,\,0\,,\quad t>0\,,
$$
and consequently $\varphi(k)\overline{\psi}(k)=0\,$ for almost all $\,k\in\mathbb R_+.\,$ Due to analyticity of $\,\hat{e}(k)\varphi(k)\,$ this can be fulfilled provided that $\,\psi(k)=0\,$ almost everywhere on $\,\mathbb R_+\,$ and therefore ${\cal L}^{\bot}=\{0\}.\,$ In this way we have established that $\,{\cal L}={\cal H}\,$ and thus operator $\,L\,$ in question is completely nonselfadjoint.

\bigskip

{\bf Corollary 2.} {\it Positivity of $\,q(x)\,$ on a certain interval guarantees for the operator $\,L=L_0+V\,$ the absence of real  {\rm(}negative{\rm)} eigenvalues so that $\,\sigma_d(L)\subset\mathbb C_+.$}

\bigskip

It was manifested in [26] that results from the papers [1] and [2] enable one to construct transformation operators implementing similarity of  $\,L_c\,$ and $\,L_0\,$ under the hypothesis that $\,V(x)\,$ possesses a finite second momentum\,:
\begin{equation}
\int_0^{\infty}(1+x^2)\,|V(x)|\,dx\en<\en\infty\,.
\label{form10}
\end{equation}
In this case stationary (i.e. not expressed in terms of perturbed and unperturbed dynamics asymptotic behavior) representations for wave operators are to be given below explicitly which clarify relationship with the eigenfunction expansion problem for operator $\,L.\,$

\medskip

Provided that operator $\,L\,$ has no spectral singularities and condition (\ref{form10}) is satisfied the transformation
$$
\Psi f\,(k)\en=\en\sqrt{\frac2{\pi}}\,\frac{k}{e(k)}\int_0^{\infty}f(x)\,y(x,k)\,\,dx
$$
maps $\,{\cal H}_c\,$ onto $\,{\cal H}\,$ and $\,\ker\Psi={\cal H}_d.\,$ Mapping $\,\Psi\!:\,{\cal H}_c\to{\cal H}\,$ is invertible and its inverse is given by the formula
$$
\Psi^{-1}g\,(x)\en=\en\sqrt{\frac2{\pi}}\int_0^{\infty}g(k)\,y(x,k)\,\frac{kdk}{e(-k)}
$$
and moreover as regards non-vanishing of the denominator in the integrand here it does not matter whether dissipative operator $\,L\,$ has any spectral singularities or not.
Indeed taking imaginary part of equation (\ref{form1}) for the Jost solution $\,e(x,k)\,$ multiplied beforehand by $\,\overline{e(x,k)}\,$ in the case when $\,k^2\in\mathbb R\,$ we get the relation
$$
{\rm Im}\,\big(\,e''_{xx}(x,k)\,\overline{e(x,k)}\big)\en=\en q(x)|e(x,k)|^2,
$$
where the right-hand-side and hence the left-hand one is integrable on half-line $\,\mathbb R_+.\,$ Integrating by parts we obtain
$$
\int_0^{\infty}{\rm Im}\,\big(\,e''_{xx}(x,k)\,\overline{e(x,k)}\big)\,dx\en=\en k\,-\,{\rm Im}\,\big(e'_x(0,k)\,\overline{e(0,k)}\big)
$$
and thus given a spectral singularitry $\,\lambda=k^2\,$ one has
$$
k\en=\en\int_0^{\infty}q(x)|e(x,k)|^2\,dx\en >\en 0\,.
$$

The composition product $\,\,e(k)\Psi\,\Omega\,\Phi^{-1}\,\,$ is calculated explicitly in [7] and thus reduced to multiplication by the expression 
$$
1\,\,+\,\,k\!\int_0^{\infty}\!e^{ikx}\,y(x,k)\,V(x)\,dx\,,
$$
which is shown in [13] to coincide with Jost function. This implies the validity of

\bigskip

{\bf Statement 5.} {\it Under the condition {\rm (\ref{form10})} direct wave operator $\,\Omega\,$ intertwining dissipative operator  $\,L=L_0+V\,$ and $\,L_0\,$ coincides with the mapping $\,\,\Psi^{-1}\Phi\!:\,{\cal H}\to{\cal H}_c.\,$ Provided that operator  $\,L\,$ has no spectral singularities the composition $\,\,\Phi^{-1}\Psi\!:\,{\cal H}_c\to{\cal H}\,\,$ gives a representation for inverse wave operator $\,\widetilde{\Omega}.$
}

\bigskip

In conclusion we emphasize that stationary representation of direct wave operator readily exhibits how continuous spectrum eigenfunctions of operator $\,L_0\,$ are transformed under the "action" \,(cf. [14]) of $\,\Omega\,$ into continuous spectrum eigenfunctions of operator $\,L=L_0+V.\,$

\bigskip
\bigskip

\begin{center}
\bf References
\end{center}

\begin{enumerate}
\item
Naimark M. A. \, Investigation of the spectrum and the expansion in eigenfunctions of a nonselfadjoint differential operator of the second order on a semi-axis // Proc. Mos. Math. Soc., 1954, V.3, P.181-270.
\item
Levin B. Ja. \, Transformations of Fourier and Laplace types by means of solutions of differential equations of second order // Dokl. Math., 1956, V.106, №2, P.187-190.
\item
Akhiezer A. I., Pomeranchuk I. Ja. \, Some problems in nuclear theory, Gostechizdat, 1950.
\item
Marchenko V. A. \, Expansion in eigenfunctions of nonselfadjoint singular second-order differential operators // Sbornik Math., 1960, V.52, №2, P.739-788.
\item
Titchmarsh E. C. \, Eigenfunction expansion associated with second-order differential equations, Clarendon Press, 1946.
\item
Kato T. \, Wave operators and similarity for some nonselfadjoint operators // Math. Ann. 1966. V.162. P.258-279.
\item
Stankevich I. V. \, On linear similarity of certain nonselfadjoint operators to selfadjoint operators and on the asymptotic behavior for $t\to\infty$ of the solution of a non-stationary Schr\"odinger equation // Sbornik Math., 1966, V.69, №2, P.161-207.
\item
Phillips R. S. \, Dissipative hyperbolic systems // Trans. Am. Math. Soc., 1957, V.86, P.109-173.
\item
Friedrichs К. О. \, Symmetric positive linear differential equations // Comm. Pure Appl. Math., 1958, V.11, P.333-418.
\item
Glazman I. M. \, Direct methods of qualitative spectral analysis of singular differential operators, Israel Prog. Scientific Transl., 1965.
\item
Murtazin Kh. Kh. \, Properties of the resolvent of a differential operator with complex coefficients // Math. Notes, 1982, V.31, №2, P.118-125.
\item
Stepin S.A. \, Nonselfadjoint singular perturbations: a model of transition from a discrete to a continuous spectrum // Russ. Math. Surv., 1995, V.50, №6, P.1311-1313.
\item
Stepin S. A. \, Dissipative Schr\"odinger operator without a singular continuous spectrum // Sbornik Math., 2004, V.195, №6, P.897-915.
\item
Berezin F. A., Shubin M. A. \, The Schr\"odinger equation, Kluwer, 1991.
\item
Akhiezer N. I., Glazman I. M. \, Theory of linear operators in Hilbert space, Dover Publ., 1993.
\item
Simon B. \, Resonances in one dimension and Fredholm determinants // J. Func. Anal., 2000, V.178, P.396-420.
\item
Levitan B. M., Sargsyan I. S. \, Sturm-Liouville and Dirac operators, Kluwer, 1991.
\item
Sz\"okefalvi-Nagy B., Foias C. \, Analyse harmonique des operateurs de l'espace de Hilbert, Akademiai Kiado, 1967.
\item
Hoffman K. \, Banach spaces of analytic functions, Prentice-Hall, 1962.
\item
Kato T. \, Perturbation theory for linear operators, Springer-Verlag, 1966.
\item
Pavlov B. S. \, Spectral analysis of a dissipative singular Schr\"odinger operator in terms of a functional model // Contemp. Problems in Math. Fundamental Direct. 1991, V.65, P.95-163.
\item
Nikolskii N. K., Khrushchev S. V. \, A functional model and some problems of the spectral theory of functions // Proc. Steklov Inst. Math., 1988, V.176, P.101-214.
\item
Yafaev D. R. \, Mathematical scattering theory, Translations of mathematical monographs, v.105, AMS, 1992.
\item
Naboko S. N. \, Functional model of perturbation theory and its applications to scattering theory // Proc. Steklov Inst. Math., 1981, V.147, P.85-116.
\item
Stepin S. A. \, Wave operators for the linearized Boltzman equation in one-speed transport theory // Sbornik Math., 2001, V.192, №1, P.141-162.
\item
Liantse V. E. \, Nonselfadjoint second order differential operator on semi-axis // Appendix to the book: Naimark M. A. Linear differential operators, Nauka, 1969, P.443-498.
\end{enumerate}

\noindent Institute of Mathematics\\ University  of Bia{\l}ystok\\
 ul.~Akademicka, 2\\ PL-15-267  Bia{\l}ystok\\ Poland\\
\end{document}